\documentclass[11pt]{article}

\usepackage{amscd,amsmath, amssymb, fancyhdr, epsfig,url}

\numberwithin{equation}{section}

\newcommand{\version}{version 1.0,\ \  19.08.2017 }

\def\eqref#1{(\ref{#1})}
\newcommand{\goth}{\mathfrak}
\newcommand{\g}{{\mathfrak g}}
\newcommand{\arrow}{{\:\longrightarrow\:}}
\newcommand{\Z}{{\Bbb Z}}
\newcommand{\C}{{\Bbb C}}

\newcommand{\R}{{\Bbb R}}
\newcommand{\Q}{{\Bbb Q}}

\def\1{\sqrt{-1}\:}
\newcommand{\restrict}[1]{{\left|_{{\phantom{|}\!\!}_{#1}}\right.}}
\newcommand{\cntrct}                
{\hspace{2pt}\raisebox{1pt}{\text{$\lrcorner$}}\hspace{2pt}}

\makeatletter
\def\x@arrow{\DOTSB\Relbar}
\def\xlongequalsignfill@{\arrowfill@\x@arrow\Relbar\x@arrow}
\newcommand{\xlongequal}[2][]{%
        \ext@arrow 0099\xlongequalsignfill@{#1}{#2}}
\def\xlongrightarrowfill@{\arrowfill@\relbar\relbar\longrightarrow}
\newcommand{\xlongrightarrow}[2][]{%
        \ext@arrow 0099\xlongrightarrowfill@{#1}{#2}}
\makeatother

\renewcommand{\bar}{\overline}
\renewcommand{\phi}{\varphi}
\renewcommand{\epsilon}{\varepsilon}
\renewcommand{\geq}{\geqslant}
\renewcommand{\leq}{\leqslant}

\newcommand{\Gr}{\operatorname{Gr}}

\newcommand{\Cl}{\operatorname{Cl}}

\newcommand{\Pos}{\operatorname{Pos}}
\newcommand{\Kah}{\operatorname{Kah}}

\renewcommand{\Re}{\operatorname{Re}}

\newcommand{\Teich}{\operatorname{\sf Teich}}

\newcommand{\Per}{\operatorname{\sf Per}}
\newcommand{\Perspace}{\operatorname{{\Bbb P}\sf er}}

\newcounter{Mycounter}[section]
\newcounter{lemma}[section]
\newcounter{claim}[section]
\renewcommand{\theclaim}{{Claim \thesection.\arabic{claim}}}
\newcommand{\claim}{%
    \setcounter{claim}{\value{Mycounter}}
    \refstepcounter{claim}
    \stepcounter{Mycounter}
    {\noindent \bf \theclaim:\ }}

\newcounter{sublemma}[section]
\newcounter{corollary}[section]
\newcounter{theorem}[section]
\renewcommand{\thetheorem}{{Theorem \thesection.\arabic{theorem}}}
\newcommand{\theorem}{%
    \setcounter{theorem}{\value{Mycounter}}
    \refstepcounter{theorem}
    \stepcounter{Mycounter}
    {\noindent \bf \thetheorem:\ }}

\newcounter{conjecture}[section]
\newcounter{proposition}[section]
\renewcommand{\theproposition}
      {{Proposition \thesection.\arabic{proposition}}}
\newcommand{\proposition}{%
    \setcounter{proposition}{\value{Mycounter}}
    \refstepcounter{proposition}
    \stepcounter{Mycounter}
    {\noindent \bf \theproposition:\ }}

\newcounter{definition}[section]
\renewcommand{\thedefinition}
      {{Definition~\thesection.\arabic{definition}}}
\newcommand{\definition}{%
    \setcounter{definition}{\value{Mycounter}}
    \refstepcounter{definition}
    \stepcounter{Mycounter}
    {\noindent \bf \thedefinition:\ }}

\newcounter{example}[section]
\newcounter{remark}[section]
\renewcommand{\theremark}{{Remark \thesection.\arabic{remark}}}
\newcommand{\remark}{%
    \setcounter{remark}{\value{Mycounter}}
    \refstepcounter{remark}
    \stepcounter{Mycounter}
    {\noindent \bf \theremark:\ }}

\newcounter{problem}[section]
\newcounter{question}[section]
\makeatletter

\@addtoreset{equation}{section} \@addtoreset{footnote}{section} \makeatother

\def\blacksquare{\hbox{\vrule width 5pt height 5pt depth 0pt}}
\def\endproof{\blacksquare}

\begin{document}
\begin{center}
{\LARGE\bf
Ergodic complex structures on hyperk\"ahler manifolds: an erratum\\[4mm]
}

Misha Verbitsky\footnote{Misha Verbitsky is partially supported by the 
Russian Academic Excellence Project '5-100'.}

\end{center}

{\small \hspace{0.10\linewidth}
\begin{minipage}[t]{0.85\linewidth}
{\bf Abstract} \\
Let $M$ be a hyperk\"ahler manifold, $\Gamma$
its mapping class group, and $\Teich$ the Teichm\"uller
space of complex structures of hyperk\"ahler type.
After we glue together birationally equivalent points,
we obtain the so-called birational Teichm\"uller space $\Teich_b$.
Every connected component of $\Teich_b$ is identified
with $\Perspace=\frac{SO(3, b_2-3)}{SO(2)\times SO(1, b_2-3)}$
by global Torelli theorem. The mapping class group of $M$
acts on $\Perspace$ as a finite index subgroup of the 
group of isometries of the integer cohomology lattice, that
is, satisfies assumptions of Ratner theorem. We prove
that there are three classes of orbits, closed, dense
and the intermediate class which corresponds to varieties
with $\dim \Re (H^{2,0}(M))\cap H^2(M,\Q)=1$.
The closure of the later orbits is a fixed point set
of an anticomplex involution of $\Perspace$. 
This fixes an error in the paper \cite{_V:Ergodic_},
where this third class of orbits was overlooked.
We explain how this affects the works based on
\cite{_V:Ergodic_}.
\end{minipage}
}

\tableofcontents


\section{Introduction}
\label{_Introduction_Section_}

\subsection{Orbits of monodromy action: an extra orbit}

This paper is essentially an erratum to \cite{_V:Ergodic_}.
For an introduction to Teichm\"uller theory, global Torelli
theorem and Ratner theory as applied to the Teichm\"uller
space, please see \cite{_V:Ergodic_}.

Let $M$ be a hyperk\"ahler manifold, and $\Teich$ the Teichm\"uller space
of complex structures of hyperk\"ahler type. Fix a connected component
$\Teich^I$ of $\Teich$ (there are finitely many of them, per \cite{_Huybrechts:finiteness_}), and let $\Gamma$
be the subgroup of mapping class group of $M$ mapping $\Teich^I$
to itself. If we glue together the points of $\Teich^I$ which are
non-separable, we obtain a Hausdorff manifold $\Teich_b$, which is
diffeomorphic to its period space
$\Perspace=\frac{SO(3, b_2-3)}{SO(1, b_2-3)\times SO(2)}$, as follows from the global
Torelli theorem (\cite{_V:Torelli_}). The mapping class group
of $M$ acts on $\Perspace$ as an arithmetic lattice, that is,
a finite index subgroup in $SO(H^2(M, \Z))$.

In \cite{_V:Ergodic_}, we classified the orbits of $\Gamma$
on $\Teich$ and proved that a general orbit of $\Gamma$ is dense
in $\Teich$, and an orbit which corresponds to a manifold of
maximal Picard rank is closed. However, one class of orbits
was overlooked: the orbits of $I\subset \Teich$ such that
$\Re(H^{2,0}(M,I))\cap H^2(M,\Q)$ has rank 1.

In the present erratum we give a complete classification of
orbits of $\Gamma$-action on $\Perspace$, and prove that
the closure of the intermediate orbit is totally real: each of these orbits is a fixed
point set of an anti-complex involution. In particular,
it does not contain positive-dimensional complex subvarieties,
and is not contained in any proper complex subvariety.

In the last section we check that the geometric applications obtained
using the older version of \cite{_V:Ergodic_} remain valid
with exception of one statement from \cite{_KLV_}
(not stated as a theorem, but implied in the abstract).

The third orbit was omitted in \cite{_V:Ergodic_}
because here we actually classified the orbits of
the group $H := SO(1, b_2-3)\times SO(2)$ on the space
$SO(3, b_2-3)/\Gamma$. However, one cannot apply
Ratner's theorem to the group $H$, because it is not
generated by unipotents. To use Ratner's theorem,
we need to consider orbits of the smaller group
$H':=SO(1, b_2-3)\subsetneq H$.

Ratner's theorem implies that a closure of an
orbit of $H'$ is an orbit of an intermediate
subgroup $S$, with $H' \subset S  \subset SO(3, b_2-3)$.
There are no proper intermediate subgroups between
$SO(1, b_2-3)\times SO(2)$ and $SO(3, b_2-3)$;
however, there is an intermediate subgroup $SO(2, b_2-3)$
between $SO(1, b_2-3)$ and $SO(3, b_2-3)$, and it
corresponds to the extra orbit.

\subsection{Preliminaries}

Throughout this paper, {\bf hyperk\"ahler manifold} means a compact
holomorphically symplectic manifold $M$ of K\"ahler type, satisfying
$\pi_1(M)=0, H^{2,0}(M)=\C$.  Its second cohomology
are equipped with a primitive, bilinear symmetric,
non-degenerate integral form $q$ of signature $(3, b_2-3)$,
called {\bf Bogomolov-Beauville-Fujiki form}, or ``BBF form''.
The corresponding
{\bf Teichm\"uller space} $\Teich$ is the space of all complex structures 
of hyperk\"ahler type up to isotopy.
We shall also speak of ``Teichm\"uller spaces'' for other geometric structures,
defined in the same way.
The {\bf monodromy group} associated with a component of $\Teich$
is the subgroup of $O(H^2(M, \Z)$. The {\bf period space} $\Perspace$ is
the Grassmanian of positive, oriented 2-planes in $H^2(M,\R)$.
The {\bf period map} $\Per:\; \Teich \arrow \Perspace$
takes a manifold $(M,I)$ to the 2-plane
$\Re H^{0,2}(M,\R)\subset H^2(M,\R)$.
``Picard rank'' of a hyperk\"ahler manifold
is $\dim_\Q (H^{1,1}(M)\cap H^2(M, \Q))$.
Global Torelli theorem (\cite{_V:Torelli_})
states that the map $\Per$ is ``the Hausdorff quotient'':
it is surjective, and satisfies $\Per(I)=\Per(I')$ 
if and only if $I$ and $I'$ are non-separable in $\Teich$.

The {\bf monodromy group} $\Gamma$ is the subgroup of $O(H^2(M,\Z))$
generated by monodromy of all Gauss-Manin local systems for
all complex deformations of $(M,I)$ over a compact base. 
This group acts on the corresponding connected
component of the Teichm\"uller space, denoted as $\Teich^I$,
and the quotient set $\Teich^I/\Gamma$ is precisely the set
of all equivalence classes of complex
structures on $M$ in the same deformation class as $I$.
In \cite{_V:Ergodic_} it was shown that $\Gamma$ is a finite index
subgroup of $O(H^2(M,\Z))$, and its action on $\Perspace$ is ergodic.

We use the notation
$(M,I)$ to put emphasis on the complex structure $I$
defined on the underlying smooth manifold $M$.

The connected component subgroup of the group $SO(p,q)$
is denoted as $SO^+(p,q)$.


\section{Subgroups of $SO(a,b)$ and Ratner theory}


\subsection{Lie subalgebras $\goth{so}(a-2,b)\subsetneq \goth{g} \subsetneq\goth{so}(a,b)$}

In this subsection, we classify connected Lie subgroups
 $G\subset SO(a,b)$ containing $SO^+(a-2,b)$.
 To classify connected subgroups
it suffices to classify Lie subalgebras 
$\goth{so}(a-2,b)\subsetneq \goth{g} \subsetneq\goth{so}(a,b)$.
This is done as follows.

\hfill

\theorem\label{_intermedia_subgrou_Theorem_}
Let $V= \R^{p+q}$ be a vector space equipped with a quadratic form
of signature $(p, q)$, $V_1 \subset V$ a positive 2-dimensional
subspace, and $V_0 = V_1 ^\bot$. The Lie algebra
$\goth{so}(V_0)$ is identified with a subalgebra
of $\goth{so}(V)$ acting trivially on $V_1$. Consider a proper
Lie subalgebra $\g \subsetneq \goth{so}(V)$ such that
$\goth{so}(V_0)\subsetneq \goth{g} \subsetneq\goth{so}(V)$.
Then  either $\goth g = \goth{so}(V_2)$, where $V_2 \subset V$ is a $p+q-1$-dimensional
subspace, $V_0\subsetneq V_2 \subsetneq V$, or 
$\goth g = \goth{so}(V_0)\oplus \goth{so}(V_1)$.

\hfill

{\bf Proof. Step 1:} We prove that
$\g=\goth{so}(V_0)\oplus \goth{so}(V_1)$ or
$\g$ is isomorphic to $\goth{so}(V_0)\oplus V_0$ as
an $\goth{so}(V_0)$-representation.
Consider the decomposition of $\goth{so}(V)$ as 
an $\goth{so}(V_0)$-representation:
\[
\goth{so}(V)=\goth{so}(V_0) \oplus (V_0\otimes V_1) \oplus \goth{so}(V_1).
\]
Since $\goth g\subset \goth{so}(V)$ is $\goth{so}(V_0)$-invariant, 
it is an $\goth{so}(V_0)$-submodule in $\goth{so}(V)$.
Since $\goth{so}(V_0) \oplus V_0\otimes V_1\subset \goth{so}(V)$ 
generates the Lie algebra $\goth{so}(V)$, one has $\goth g \not \supset V_0\otimes V_1$.
Then either  $\goth g \cap V_0 \otimes V_1=0$, in which case
$\goth g = \goth{so}(V_0)\oplus \goth{so}(V_1)$ (the space $\goth{so}(V_1)$
is 1-dimensional), or $\goth g \cap V_0 \otimes V_1=0$
is a proper $\goth{so}(V_0)$-invariant 
subspace of $ V_0 \otimes V_1$. Since $V_1$ is a trivial 2-dimensional
representation of the Lie algebra $\goth{so}(V_0)$, the space  $V_0 \otimes V_1$
is isomorphic (as an $\goth{so}(V_0)$-module) to a direct sum of two
copies of $V_0$, which is an irreducible representation
of $\goth{so}(V_0)$. Then any proper subrepresentation of
$V_0 \otimes V_1$ is isomorphic to $V_0$. 

{\bf Step 2:} Let $\g=\goth{so}(V_0)\oplus R$,
where $R\subset V_0 \otimes V_1$ is a $\goth{so}(V_0)$-invariant
subspace. Then $R= V_0 \otimes R_1$, where $R_1 \subset V_1$
is a 1-dimensional subspace in a 2-dimensional subspace $V_1$.
Denote by $S\subset V_1$ the orthogonal complement to $S$ in $V_1$, 
and let $V_2=V_0 \oplus S$.
Then $\goth g$ fixes $S$, which gives a Lie algebra homomorphism
$\g \arrow \goth{so}(V_2)$. Comparing dimensions, we find
that it is an isomorphism.
\endproof

\subsection{Ratner theory for $SO^+(a-2,b)\subset SO^+(a,b)$}

To fix the notation, 
we state the Ratner orbit closure theorem in a particular situation
where we needed it. For a general version, see \cite{_V:Ergodic_} 
or \cite{_Morris:Ratner_}. As usual, $SO^+(a,b)$ denotes the
connected component of the special orthogonal group $SO(a,b)$.

\hfill

\theorem\label{_Ratner_for_SO(a,b)_Theorem_}
Consider an integer lattice $V_\Z=\Z^{a+b}$ equipped with an integer
scalar product $q$ of signature $(a,b)$, $a>2$, $b>0$, $a+b >4$. Let $\Gamma\subset SO(V_\Z)$
be a finite index sublattice, $V_\R:= V_\Z\otimes_\Z \R$, and $SO^+(V_\R, q)$
the connected component of its orthogonal group. Replacing $\Gamma$
by a finite index subgroup if necessarily, we may assume
that $\Gamma$ belongs to the connected component $SO^+(V_\R)$
of the group $SO(V_\R)$. Let $H\subset SO^+(V_\R, q)$
be the group of all elements of $SO^+(V_\R)$ acting trivially  on a
2-dimensional positive subspace $V_1 \subset V_\R$.
Consider the left action of $H$ on $G/\Gamma$; let $Hx$
be an orbit of $x\in G/\Gamma$, and $\overline{Hx}$
its closure in standard topology. Then there exists
a rational point $y\in \overline{Hx}$ and a rational Lie subgroup
$H_1$ generated by unipotents, $H\subset H_1 \subset G$, such that $\overline{Hx}= H_1 y$.

\hfill

{\bf Proof:}
This is the statement of Ratner theorem from  \cite{_Morris:Ratner_};
to apply it, we have to check that the group $H$ is generated by unipotents.
The universal cover of $SO^+(2,2)$  is $SL(2,\R)\times SL(2, \R)$, hence
it is generated by unipotents. For other values of $(a,b)$, $a+b > 4, a\geq 3, b > 1$,
the group $H=SO^+(a-2, b)$ is simple, and its subgroup generated by
unipotents is normal, hence it is equal to $H$ or empty. There are
unipotent elements in $SO^+(a,b)$ for $a\geq 2, b \geq 1$,
which gives the statement of \ref{_Ratner_for_SO(a,b)_Theorem_}. \endproof

\hfill

\remark
In \cite[Example 4.4]{_V:Ergodic_} the Ratner Theorem was stated in the same
way and then applied to $H=SO(a-2,b)\times SO(2)$.
However the group $SO(a-2,b)\times SO(2)$ is not generated by
unipotents. Its smallest subgroup generated by unipotents
is $SO(a-2,b)$. This is why in place of \cite[Example 4.4]{_V:Ergodic_}
we should use \ref{_Ratner_for_SO(a,b)_Theorem_}.

\subsection{Ratner theory and complex geometry of the period space}

\claim \label{_Perspace_complex_Claim_}
Let $V= \R^{a,b}$ be a vector space equipped with a non-degenerate scalar
product $q$ of signature $(a,b)$, and $\Gr_{++}= \frac{SO^+(a,b)}{SO^+(a-2, b)\times SO(2)}$
the Grassmannian of positive, oriented 2-planes. Then $Gr_{++}$ is in bijective correspondence with
\[
\Perspace:= \{l\in {\Bbb P}V_\C\ \ | \ \  q(l,l)=0, q(l, \bar l) >0\}.
\]
where $V_\C$ is the complexification of $V$.

{\bf Proof:} See, for example, \cite[Proposition 2.12]{_V:Ergodic_}. \endproof

\hfill

\theorem\label{_closure_orbit_Gr++_Theorem_}
Consider an integral lattice $V_\Z=\Z^{a+b}$ equipped with an integer-valued
scalar product $q$ of signature $(a,b)$, $a>2$, $b>0$, $a+b >4$, $V:= V_\Z\otimes_\Z \R$,
$\Gamma\subset SO(V_\Z)$ a finite index sublattice, $l\in \Gr_{++} (V)$ a point, 
and $\Gamma l\subset \Gr_{++} (V)$ its orbit. Then one of the following three possibilities is true.
\begin{description}
\item[(i)] $\Gamma l$ is closed. This happens when the 2-plane $l\subset V$
is rational. 
\item[(ii)] $\Gamma l$ is dense in $\Gr_{++} (V)$. This happens when $l$ contains
no rational vectors.
\item[(iii)] The closure $\overline{\Gamma l}$  is the set of all 2-planes
$V_1\in \Gr_{++}(V)$ containing a 1-dimensional rational subspace $v\subset V_\Q$.
\end{description}
{\bf Proof:}
Let $P:= SO(V)/SO(V_0)$, where $V_0\subset V$ is a codimension 2 subspace
of signature $(a-2,b)$. Clearly, $P$ is fibered over $\Gr_{++}(V)$
with the fiber $SO(2)=S^1$:\[
\pi:\; P \xlongrightarrow{\ /SO(2)\ }{} \frac{SO^+(a,b)}{SO^+(a-2, b)\times SO(2)}=\frac{P}{SO(2)}.
\] Since $SO(2)$ is compact, the closure of
a $\Gamma$-orbit of $x\in P$ satisfies 
\begin{equation}\label{_closure_orbits_SO(2)_Equation_}
\pi(\overline{\Gamma\cdot l})=\overline{\Gamma\cdot \pi(l)}.
\end{equation}

Notice that an $H$-invariant subset $Z$ on $G/\Gamma$ 
is closed if and only if its image is closed in the
double quotient $H\backslash G/\Gamma$. Therefore, 
the closure of an orbit $H\cdot x$ in $G/\Gamma$ 
can be obtained by taking the closure of 
$x\cdot \Gamma$ in the left quotient $H\backslash G$.
Using Ratner's theorem (\ref{_Ratner_for_SO(a,b)_Theorem_}),
we find that such closures correspond to intermediate
subgroups  $H \subset S \subset G$. For $\Gamma$,
$G$ and $H$ as in \ref{_Ratner_for_SO(a,b)_Theorem_},
the list of possible $S$ is provided by \ref{_intermedia_subgrou_Theorem_}.

To obtain a classification of $\Gamma$-orbits in $\Gr_{++}$,
we consider the $\Gamma$-action on $\frac{SO^+(a,b)}{SO^+(a-2, b)}$,
and use \ref{_intermedia_subgrou_Theorem_},
\eqref{_closure_orbits_SO(2)_Equation_} and
\ref{_Ratner_for_SO(a,b)_Theorem_}. From 
\eqref{_closure_orbits_SO(2)_Equation_} it is clear
that it suffices to classify the orbits of $\Gamma$-action
on $\frac{SO^+(a,b)}{SO^+(a-2, b)}$. From \ref{_Ratner_for_SO(a,b)_Theorem_},
it follows that the classes
of orbits in \ref{_Ratner_for_SO(a,b)_Theorem_} 
correspond to classes of intermediate
rational subgroups $S\subset SO^+(a,b)$ generated by unipotents,
with $SO^+(a-2, b)\subset H\subset SO^+(a,b)$.
From \ref{_intermedia_subgrou_Theorem_}, we obtain that
$S$ is either $SO^+(a-2, b)$, $SO^+(a-1, b)$, $SO^+(a, b)$ or
$SO^+(a-2, b)\times SO(2)$; however, the later case is 
impossible, because $S$ is generated by unipotents.
The remaining three classes give 
the three cases (i)-(iii) of \ref{_closure_orbit_Gr++_Theorem_}. 
\endproof

\hfill

The last class of orbits (\ref{_closure_orbit_Gr++_Theorem_}, (iii))
is what is missing from \cite{_V:Ergodic_}. However, most applications
of \cite{_V:Ergodic_}, obtained since, remain valid, because the extra
orbit has sufficiently high codimension. This is explained in more detail
in Section \ref{_appli_Section_}. Also, the extra orbits are negligible from the
complex-geometric point of view, which is implied by the following
result. Recall that a subvariety $N\subset M$ of a complex manifold
is called {\bf totally real} if there exists an anti-holomorphic involution
$\tau:\; M \arrow M$ such that $N$ is its fixed point set. Totally real 
subvarieties of smooth varieties are always
smooth (the fixed point set of an involution is always smooth).
From the definition it is easy to see that $TM\restrict N =TN \otimes_\R \C$.
Therefore, a tangent 
space $T_xN$ to a totally real submanifold contains no positive-dimensional
complex subspaces, and it is not contained in a proper complex subspace of $T_xM$.

This gives the following claim.

\hfill

\claim
Let $M$ be a complex manifold, and $N\subset M$
a totally real submanifold. Then $N$ contains no positive-dimensional
complex subvarieties of $M$, and is contained in no proper
subvariety of $M$.
\endproof

\hfill

\proposition\label{_totally_real_Gr(v)_Proposition_}
Let $V$ be a vector space equipped with a scalar product of signature
$(a, b)$, $a>2$, and $\Gr_{++}:=\frac{SO^+(a,b)}{SO^+(a-2, b)\times SO(2)}$
the positive oriented Grassmannian, equipped with the complex structure
as in \ref{_Perspace_complex_Claim_}. Fix a vector $v\in V$,=
with $(v,v)>0$, and let $\Gr_{++}(v)\subset \Gr_{++}$
be the set of all positive, oriented  2-planes containing $v$.
Then $\Gr_{++}(v)$ is a fixed point set of an antiholomorphic
involution of $\Gr_{++}$. This means that $\Gr_{++}(v)$
is totally real, contains no positive-dimensional
complex subvarieties of $\Gr_{++}$, and is contained in no proper
complex subvariety of $\Gr_{++}$.

\hfill

{\bf Proof:}
Let $\iota_v:\; V \arrow V$ be the reflection
\[
x \mapsto x- 2\frac{(x,v)}{(v,v)} v.,
\]
and $\gamma_v:\; \Gr_{++}\arrow \Gr_{++}$
be the composition of the map induced by $\iota_v$
and the orientation reversal. Clearly, $\gamma_v$ fixes
a 2-plane $W\subset V$ if and only if $W$ contains $v$.
Since the complex structure on $\Gr_{++}$ is replaced
by its opposite when the orientation changes,
$\gamma_v$ is an anti-holomorphic map.
\endproof


\section{Mapping  class group action on the Teichm\"uller space}
\label{_map_Teich_Section_}

\subsection{Mapping class group action and non-Hausdorff points}
\label{_non_haus_mapping_class_Subsection_}

The period space $\Perspace$ is obtained by taking the 
Teichm\"uller space $\Teich$
and gluing together all non-separable points. It is
not obvious how to relate closed or dense orbits
on $\Perspace$ to that on $\Teich$. The answer to
this problem, given in \cite{_V:Ergodic_}, remains 
valid. Let us remind how it is done.

\hfill

\theorem\label{_Teichm_ergo_NH_Theorem_}
Let $\Per:\; \Teich\arrow \Perspace$ be the period  map,
and $\Gamma$ the mapping class group. Assume that $b_2(M)\geq 5$.
Then the period map commutes with taking closures:
$\Per(\overline{\Gamma I})=\overline{\Gamma \Per(I)}$,
for all $I\in \Teich$.

\hfill

The proof of \ref{_Teichm_ergo_NH_Theorem_}
takes the rest of this section.

Recall that {\bf the positive cone}
$\Pos([I])$ is the connected component of
the set of all real (1,1)-classes $v\in H^{1,1}([I])$
satisfying $q(v,v)>0$, where $q$ is the Bogomolov-Beauville-Fujiki
form. There are two connected components,
and we choose one which contains the K\"ahler classes.
A subset $K\subset \Pos([I])$
is called {\bf K\"ahler chamber} if it is a K\"ahler cone
for some $I\in \Teich$ satisfying $\Per(I)=[I]$.
The following fundamental result is due to Eyal Markman.

\hfill

\proposition\label{_chambers_periods_teich_Proposition_}
Different K\"ahler chambers of $[I]$ do not intersect, and
$\Pos([I])$ is a closure of their union. Moreover, there is
a bijective correspondence between points of
$\Per^{-1}([I])$ in the Teichm\"uller component of $I$
and the set of K\"ahler chambers of $[I]$.

{\bf Proof:} \cite[Proposition 5.14]{_Markman:survey_}. \endproof

\hfill

Consider the set $\Teich_K$ of pairs $\left(I\in \Teich, \omega\in \Kah(M,I)\right)$,
where $\Kah(M,I)$ denotes the K\"ahler cone, and $q(\omega, \omega)=1$.
 Let $\Per_K$ be the set of
all pairs 
\[ \{\left([I]\in \Perspace, \omega\in \Pos([I])\right) \ \ |\ \ q(\omega, \omega)=1\}.
\]
Consider the period map $\Per_K:\; \Teich_K\arrow \Perspace_K$ mapping
$(I,\omega)$ to $(\Per(I),\omega)$. By \ref{_chambers_periods_teich_Proposition_},
$\Per_K$ is injective with dense image. This means that points
of $\Teich$ can be identified with the space ${\goth P}$ pairs $(I, k)$,
where $I\in \Perspace$ and $k\subset H^{1,1}(M,I)$ is 
one of the K\"ahler chambers of $I\in \Perspace$. The
topology on this space of pairs is the quotient topology
induced by the map $\tau:\; \Perspace_K \arrow {\goth P}$.

Taking in account this picture, the closure of the monodromy orbit of $I\in \Teich$ 
can be understood as the closure of $\Gamma\cdot (\Per(I),k)$ 
in $\Teich_K\subset \Perspace_K $, where
$k$ is the K\"ahler chamber $I\in \Teich$, and $\Per(I),k$
is considered as a subset in $\Perspace_K$.

To prove that $\overline{\Gamma \cdot (\Per(I),k)}$
contains the positive cone $(I, \Pos(I))\subset \Perspace_K$,
we use the Ratner's theorem again.
We find an orbit $R$ of a one-parametric group $H_0\subset SO(H^2(M,\R))$ generated by unipotents
and fixing $I\in \Perspace$
with $R\subset \Perspace_K$ lying within $(\Per(I),k)$, and show
(using Ratner's theorem) that $\Gamma \cdot R$ is dense in $\Perspace_K$.

The same strategy was used in \cite{_V:Ergodic_},
however, the one-parametric group that we chose in \cite{_V:Ergodic_}
was in fact not generated
by unipotents. Here we fill this gap and give a corrected version
of the proof.

\subsection{Round segments on the boundary of the K\"ahler cone}

In this subsection we shall speak of ``round bits'' on the boundary
of the convex cone. The geometric intuition behind this terminology
is that the K\"ahler cone is polyhedral at some points, and strictly
convex with smooth boundary at other points; the points of the boundary
where it's smooth and strictly convex constitute ``round bits''.

\hfill

\definition
Let $(M,I)$ be a hyperk\"ahler manifold, and $W\subset H^{1,1}_I(M,\R)$
a 3-dimensional real subspace of signatue $(1,2)$. We say that
{\bf the K\"ahler cone boundary has round bits on $W$}
if the projectivization of the
closure ${\Bbb P}\overline{\Kah(M,I)\cap W}$ contains an
open subset of the boundary of the disk
${\Bbb P}\Pos(W)$.\footnote{Here, as elsewhere,
  $\Pos(W)$ denotes the positive cone,
  which is one of two components of the set $\{w\in W \ \ |\ \ q(w,w)>0\}$.
  Notice that ${\Bbb P}\Pos(W)$ is
  a 2-dimensional disk in ${\Bbb P}(W)=\R P^2$, identified
  with one of the standard models of the hyperbolic plane.}

\hfill

The following proposition
is used to construct a one-parametric unipotent orbit which lies
in $\Kah(M,I)$. The idea is to construct a horocycle in the hyperbolic segment
${\Bbb P}\Pos(W)\subset {\Bbb P}\Kah(M,I)= (I,k)\subset \Perspace_K$.
However, a horocycle exists only in subsets of $\overline{\Pos(W)}$
which contain open parts of the boundary.

\hfill

\proposition\label{_rnd_bits_exists_Proposition_}
Let $(M,I)$ be a hyperk\"ahler manifold such that
$b_2(M)\geq 5$ and the Picard rank of $(M,I)$ is non-maximal.
Denote by ${\goth T}$ the orthogonal complement to
$H^{1,1}_I(M,\R)\cap H^2(M,\Q)$ in $H^{1,1}_I(M,\R)$;
since the Picard rank of $(M,I)$ is non-maximal, the space
${\goth T}$ has positive dimension.
Then for each 3-dimensional space $W\subset H^{1,1}_I(M,\R)$
of signature $(1,2)$ such that ${\goth T}\cap W\neq 0$,
the K\"ahler cone boundary has round bits on $W$.

\hfill

{\bf Proof:} By \cite[Theorem 1.19]{_AV:MBM_}, 
the K\"ahler cone $\Kah(M,I)$ is a connected component of
$\Pos(M,I)\backslash S^\bot$, where $S^\bot = \bigcup_i s_i^\bot$
is a countable union of hyperplanes obtained as orthogonal complement
to the so-called MBM classes $s_i \in H^{1,1}_I(M,\Z)$.
Clearly, all $s_i$ are orthogonal to ${\goth T}$.
Therefore, all hyperplanes $s_i^\bot$ intersect
the subspace $P:= {\goth T}\cap W\subset W$. 
The corresponding partition of the disk ${\Bbb P}\Pos(W)$
onto pieces by $S^\bot$
depends on the signature of $P$.

We may assume that $P$
is one-dimensional; if it's 2-dimensional, one has $S^\bot_i\cap W=P$,
and the disk is partitioned onto at most two pieces, hence
each piece contains a round part of the boundary.

Notice that the set $s_i^\bot$ of walls of the K\"ahler chambers
is locally finite (\cite[Proposition 10]{_Hassett_Tschinkel:moving_}).
On the disk ${\Bbb P}\Pos(W)$ these lines trace geodesics which
intersect in one point of the disk if ${\goth T}$ has positive
signature, and intersect in ``imaginary'' (negative) point
outside of the disk if the signature is negative,
We obtain that the hyperplanes $s_i^\bot$ partition the disk
${\Bbb P}\Pos(W)$ as follows.

\centerline{\begin{tabular}{cc}
    \includegraphics[width=0.30\textwidth]{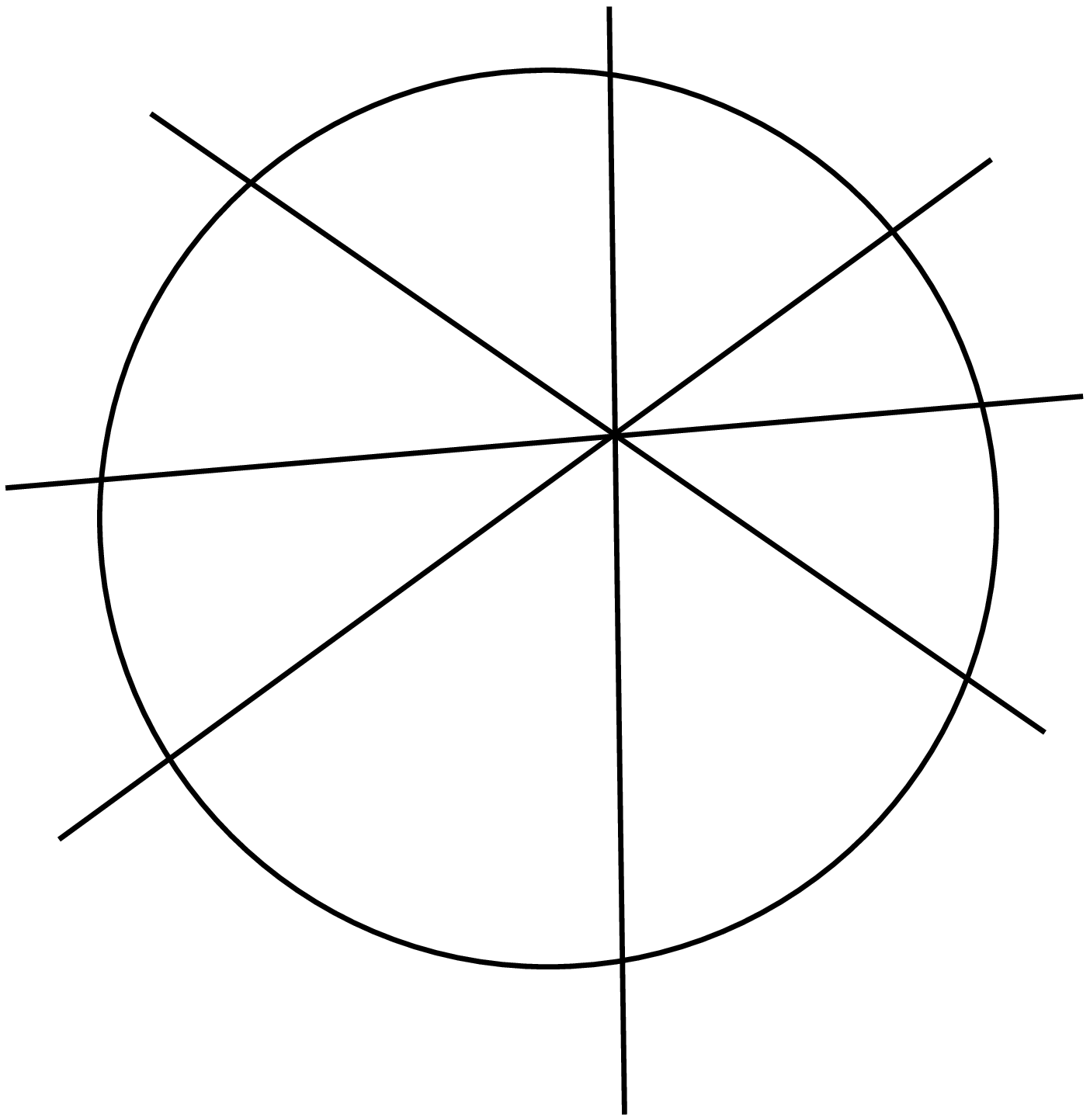}\  &\
    \includegraphics[width=0.30\textwidth]{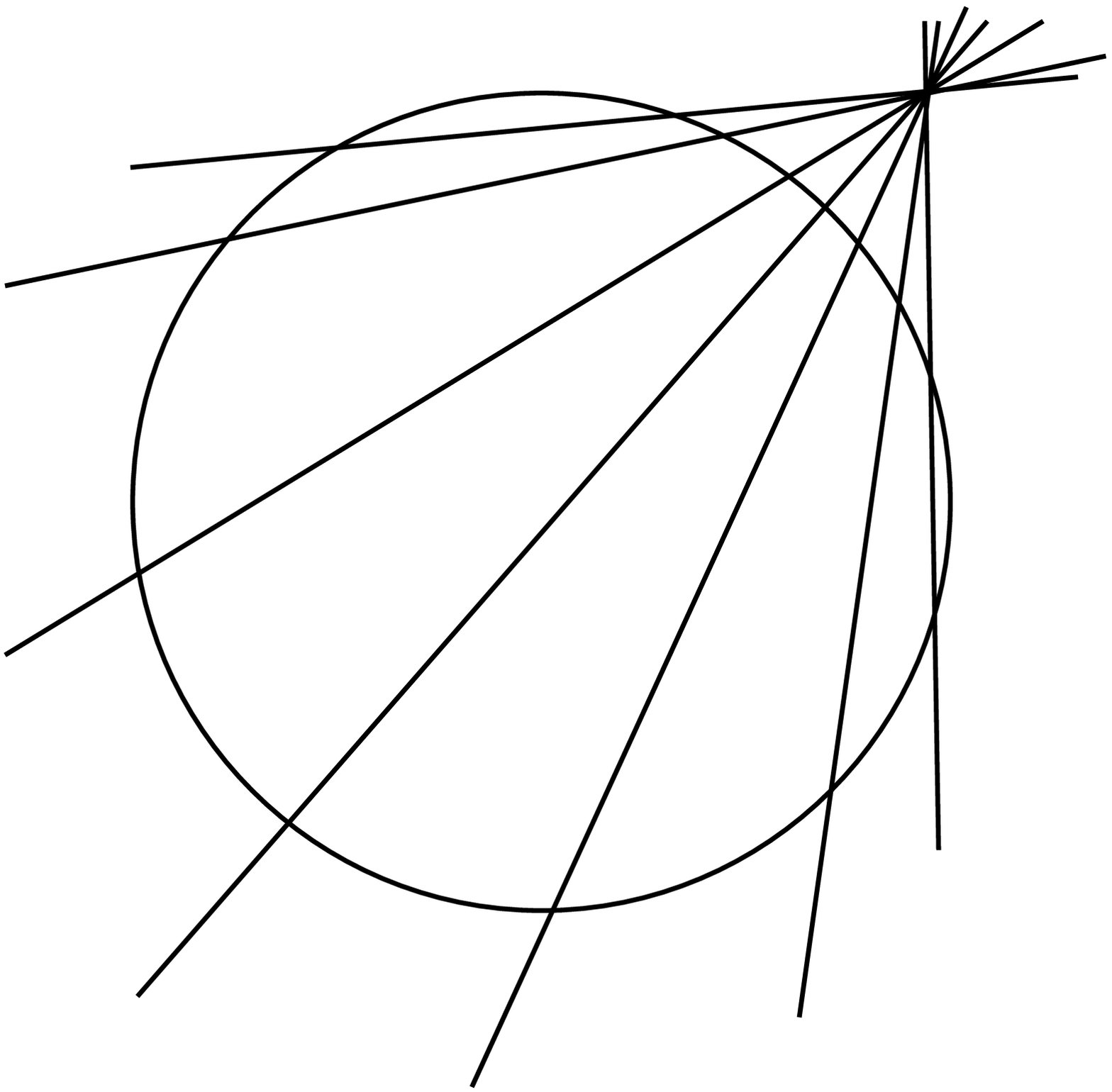} \\
{\bf \scriptsize $\goth T$ has positive signature} & {\bf \scriptsize $\goth T$ has  negative signature} 
\end{tabular}}

In both cases, as seen from the picture, all connected components
of ${\Bbb P}\Pos(W)\backslash S^\bot$ have round bits on the boundary of the disk.
\endproof

\subsection{Horocycles in the K\"ahler cone}

Let $\Perspace_K$  the period
space of K\"ahler classes (Subsection \ref{_non_haus_mapping_class_Subsection_}),
and $\Gamma$ the monodromy group acting on $\Perspace_K$.
We find a one-parametric group $H_0\subset SO^+(H^2(M,\R))$ generated by unipotents
such that for some $x\in \Kah(M,I)$,
the orbit $H_0 \cdot x$ is contained in $\Kah(M,I)$.
We apply \ref{_rnd_bits_exists_Proposition_} and find $W\subset H^{1,1}(M,\R)$
such that the K\"ahler cone boundary of $(M,I)$ has round bits on $W$.
The group $H_0$ acts on the Poincare disk $\Delta:={\Bbb P}\Pos(W)$ by
parabolic isometries, and its orbits are horocycles. From the picture below
it is clear that there are uncountably many horocycle
subgroups with orbits in $\Kah(M,I) \cap W$ whenever
the K\"ahler cone boundary of $(M,I)$ has round bits on $W$.
\begin{center}
    \includegraphics[width=0.65\textwidth]{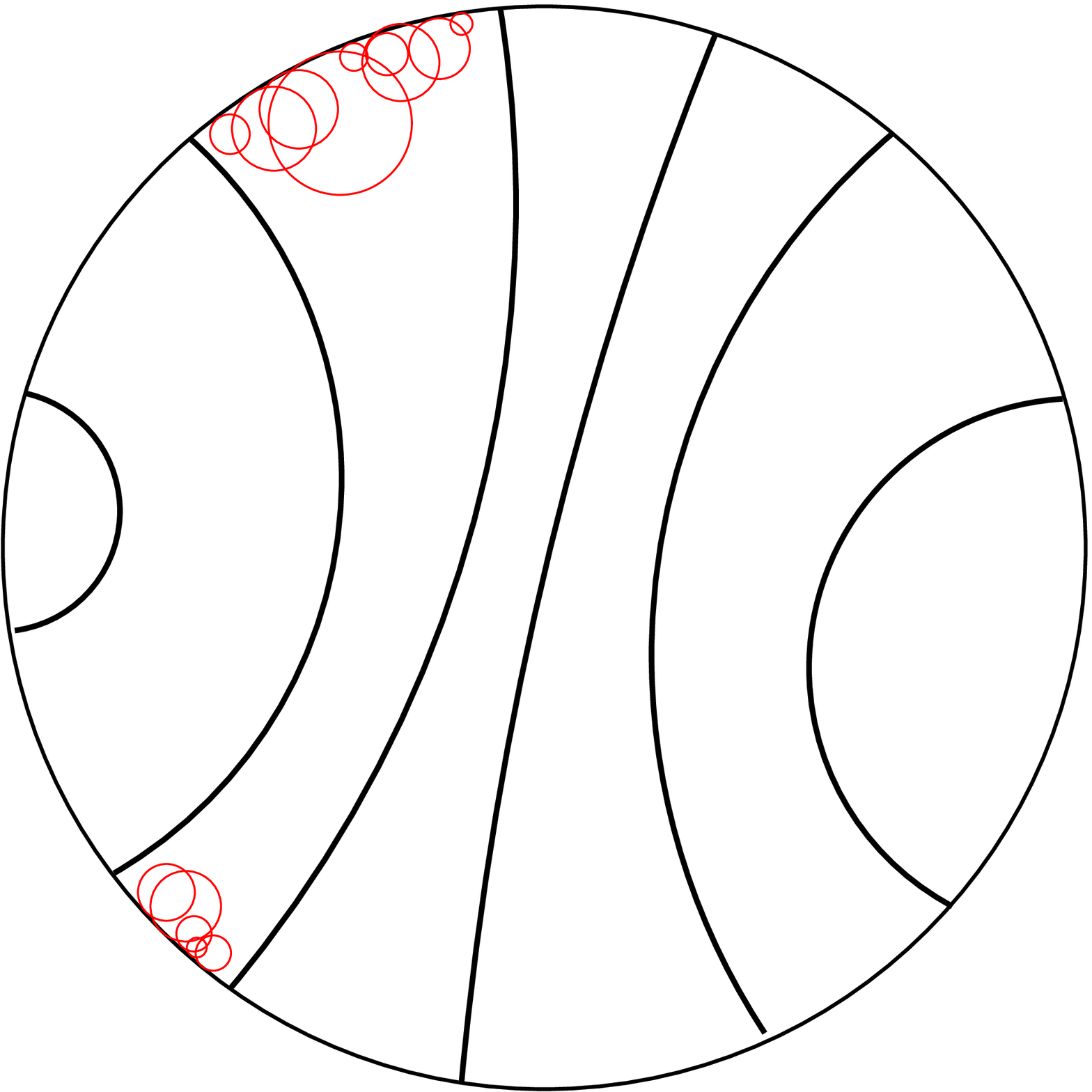}\\
{\bf \scriptsize Horocycles on a Poincare disc} 
\end{center}
The following proposition, applied together with Ratner's theorem, finishes the proof of 
\ref{_Teichm_ergo_NH_Theorem_}, as indicated in
Subsection \ref{_non_haus_mapping_class_Subsection_}.

\hfill

By ``general'' we always mean ``not contained in a union of countable
number of Zariski closed subsets of positive codimension''.

\hfill

\proposition\label{_minim_cont_H_0_Proposition_}
Let $V_\Q$ be a rational vector space of signature
$(3,d)$, $d\geq 2$, $V:=V_\Q \otimes_\Q \R$,
and $\Gamma\subset SO^+(V)$ an arithmetic lattice. Consider a subspace $V_1\subset V$
of signature $(1,d)$, let $P\subset V_1$ be a 1-dimensional subspace, $W\supset P$
a general 3-dimensional subspace of $V_1$ of signature $(1,2)$ and containing $P$, and
$H_0\subset SO(W) \subset SO(V)$
a general one-parametric unipotent subgroup. Then the smallest rational Lie subgroup
containing $H_0$ also contains $SO^+(V_1)$.

\hfill

{\bf Proof:} Denote by $\Cl(H_0)$ the smallest rational Lie subgroup
containing $H_0$, and let ${\goth I}_{H_0}\subset \coprod {\Bbb P} V^{\otimes i}$
be the set of all rational projective invariants of $H_0$. By Chevalley's theorem
(see e.g. \cite{_Morris:Ratner_}), $\Cl(H_0)$
is the biggest Lie group fixing all points $x\in{\goth I}_{H_0}$.
Since the set $ {\goth I}_{H_0}$ is countable and semicontinuous in $H_0$,
the group $\Cl(H_0)$ for generic $H_0$ is independent on
$H_0$ and contains all unipotent sugroups
of $SO^+(W)$, for all $W$ of appropriate signature satisfying
$P\subset W \subset V_1$. This implies $\Cl(H_0)\supset
SO^+(V_1)$. \endproof

\hfill

From \ref{_minim_cont_H_0_Proposition_},
\ref{_Teichm_ergo_NH_Theorem_} follows immediately.
Indeed, apply  \ref{_minim_cont_H_0_Proposition_}
to $V=H^2(M, \R)$, $V_1=H^{1,1}_I(M,\R)$, and
$P\subset{\goth T}$, and choose generic
$H_0\subset SO(W)$ which has an orbit in
$\Kah(M,I)$. By Ratner's theorem and
\ref{_minim_cont_H_0_Proposition_},
the closure of this
orbit in $\Perspace_K/\Gamma$ contains
an orbit of $SO^+(V_1)$ which contains
the set of unit vectors in the positive cone of $I$
considered as a subset of $\Perspace_K$. However,
the closure of this set in $\Perspace_K/\Gamma$
is mapped surjectively to the closure of $I$
in $\Perspace/\Gamma$, because
$\Perspace_K$ is fibered over $\Perspace$
with fiber $\Pos(M,I)$, and this fibration
is $\Gamma$-equivariant.


\section{Applications to hyperk\"ahler geometry: Hyperbolicity and Kobayashi pseudometric}
\label{_appli_Section_}

Let $\Gamma$ be the mapping class group of a complex manifold.
Recall that a complex structure is called {\bf ergodic}
if its $\Gamma$-orbit is dense in its connected component
of the Teichm\"uller space of complex structures.

In \cite{_V:Ergodic_}, we proved that a hyperk\"ahler
manifold is not Kobayashi hyperbolic. To prove this, one
needs to produce a deformation which is not Kobayashi
hyperbolic and has ergodic complex structure.

As shown by F. Campana in \cite{_Campana:twistors_},
for any hyperk\"ahler structure, the corresponding
twistor deformation has at least one non-hyperbolic fiber.
The twistor deformation can be interpreted in terms of
the period space as follows. Let $W\subset H^2(M,\R)$
be the 3-dimensional space generated by three
K\"ahler forms $\omega_I, \omega_J, \omega_K$,
associated with the quaternionic triple of complex
structures. Consider the space $S_W$ of oriented
2-planes in $W$. We can consider $S_W$ as a
subset in $\Gr_{++}$. This subset is precisely
the set of all complex structures in the twistor
deformation associated with $I,J,K$.

Let now $W_0\subset W$ be the 2-plane
corresponding to the non-hyperbolic complex structure 
constructed in \cite{_Campana:twistors_}.
If we chose $W$ such that it contains
no rational vectors, $W_0\cap H^2(M,\Q)=0$,
hence the associated $\Gamma$-orbit is dense
in $\Perspace$ (\ref{_closure_orbit_Gr++_Theorem_}).
The corresponding complex structure is ergodic.
This proves the non-hyperbolicity of hyperk\"ahler manifolds.

In \cite[Theorem 1.1]{_KLV_}, it was shown that a hyperk\"ahler
manifold $M$ with $b_2(M)\geq 7$ 
has vanishing Kobayashi pseudometric if the ``SYZ
conjecture'' holds. SYZ conjecture, which is formulated 
further on in this paper, is
the same as Kawamata's abundance conjecture
for hyperk\"ahler manifolds. It is known for the Hilbert
scheme of K3 surfaces and all its deformations
(\cite{_Bayer_Macri_}) and for the generalized Kummer
varieties (\cite{_Yoshioka_}).

By a theorem of Matsushita (see \cite{_Matsushita:fibred_}), any holomorphic
map $\pi:\; M \arrow B$ is a Lagrangian fibration
if $0 < \dim B < \dim M$. The space $B$
is known to be a Fano manifold, and it has
the same rational cohomology as $\C P^n$ if it is normal
(\cite{_Matsushita:CP^n_}).
The pullback $\eta:= \pi^*\omega$ of an ample class 
$\omega\in H^2(B,\Z)$ is an integer class on
the boundary of the K\"ahler cone of $M$
which satisfies $\int_M \eta^{\dim_\C M}=0$.
Such cohomology classes are known as
{\bf parabolic}\footnote{In more generality, parabolic classes
can be defined as nef classes of volume 0.}, and SYZ conjecture claims
that, conversely, any parabolic class 
$\eta\in H^2(M, \Z)$ is obtained this way.

Consider a manifold $M$ which has two transversal
Lagrangian fibrations (possibly rational).
Since the fibers of these fibrations are tori,
the Kobayashi pseudometric on $M$ vanishes.

To prove that any manifold $(M,I)$ 
has vanishing Kobayashi pseudometric, we need to show that
in the closure $\overline{\Gamma I}$ of its $\Gamma$-orbit there exists a
manifold with vanishing Kobayashi pseudometric.
Since the diameter of the Kobayashi pseudometric
is upper semicontinuous (\cite{_KLV_}),
this would imply that diameter of Kobayashi pseudometric
on $(M,I)$ also vanishes.

Meyer's theorem states that any integer lattice
which is not sign-definite and of rank $\geq 5$
represents 0, that is, contains an integer (1,1)-class with square 0.
Applying this result and global Torelli theorem to a
hyperk\"ahler manifold  with $b_2(M)\geq 5$, 
we can find a complex structure with any
Picard rank $\leq b_2(M)-2$ and a vector
with square 0, which would lie in a boundary
of the K\"ahler cone for some complex structure
(\cite{_AV:cusps_cone_}), and give a Lagrangian fibration
if SYZ conjecture is true.

Using the same argument for two non-proportional nef vectors, we would get
two transversal Lagrangian fibrations. If we want to have
a point with the two Lagrangian fibrations on a
positive-dimensional orbit $\Gr_{++}(v)$
(\ref{_totally_real_Gr(v)_Proposition_}), we would need to
find two Lagrangian fibration (that is, two nef classes)
on a manifold with $\Re(H^{2,0}) \ni v$,
which is done using the same argument.
This argument takes care of \cite[Theorem 1.1]{_KLV_}.

A special case of this result,
\cite[Corollary 2.2]{_KLV_}, says that any non-algebraic
K3 surface has vanishing Kobayashi metrics (for algebraic
K3 it was already known). The idea is to show that
an algebraic complex structure can be found in the 
closure of the mapping class orbit of a non-algebraic
complex structure and use the semicontinuity of
the diameter of Kobayashi metric. This argument remains true.
However, we should take care about an extra orbit
and notice that the closure of this orbit
is fixed point set of a rational anticomplex involution
on the period space. It is easy to see that this fixed point set 
intersects the divisor $\Teich_g\subset \Teich$ representing polarized
K3 of any given genus. This implies that
there is a dense set of algebraic complex
structures in any intermediate orbit.

The rest of results of \cite{_KLV_} remain the same, with
the same proofs.

\hfill

{\bf Acknowledgements:} The gap in the classification of orbits
in \cite{_V:Ergodic_} was found when I collaborated with Nessim Sibony on a different
project extending the classification of
\cite{_V:Ergodic_} to the space of complex structures equipped
with a parabolic nef class. I am grateful to Nessim for his questions
and many interesting discussions. Apologies to those who were misled 
by the omission of the elusive orbit. Many thanks to Ekaterina Amerik
and Ljudmila Kamenova for their remarks and questions about the
preliminary version of this paper. Section \ref{_map_Teich_Section_}
is based on a long series of talks with Ekaterina Amerik, who is
responsible for finding many bugs in earier versions, and many
ideas of the proof and its presentation.

\hfill

{\small

  \noindent  Misha Verbitsky\\
             {\sc Instituto Nacional de Matem\'atica Pura e Aplicada\\ Estrada Dona Castorina, 110\\
Jardim Bot\^anico, CEP 22460-320\\
Rio de Janeiro, RJ - Brasil}\\ also: \\
{\sc 
Laboratory of Algebraic Geometry, \\
Faculty of Mathematics, National Research University HSE,\\
7 Vavilova Str. Moscow, Russia,\\ 
\tt verbit@verbit.ru, verbit@impa.br }
}

\end{document}